\documentclass[11pt]{article}

\usepackage{fullpage}
\usepackage{authblk}
\usepackage{natbib}
\usepackage{algorithm}
\usepackage{algpseudocode}
\usepackage{amsmath,amsthm,amsfonts}
\usepackage{amsmath}
\usepackage[subnum]{cases}
\usepackage{tcolorbox}
\usepackage{mathrsfs}
\usepackage{amssymb}
\usepackage{color}
\usepackage{mathrsfs}
\usepackage{enumitem}
\usepackage{bm}
\usepackage{multirow}
\usepackage{makecell}
\usepackage{graphicx}
\usepackage{subcaption}
\usepackage{comment}
\usepackage{cases}
\usepackage{appendix}
\usepackage{tikz}
\usetikzlibrary{arrows,shapes}
\usepackage[colorlinks= true, linkcolor=red, citecolor=blue, urlcolor=black]{hyperref}
\usepackage{cleveref}
\usepackage{marginnote}
\usepackage{diagbox} 
\usepackage{empheq}

\numberwithin{equation}{section}

\newtheorem{theorem}{Theorem}[section]
\newtheorem{lemma}[theorem]{Lemma}
\newtheorem{corollary}[theorem]{Corollary}

\newtheorem{defn}[theorem]{Definition}

\newcommand{\mr}{\mathbb{R}}

\usepackage{booktabs}
\usepackage{pgfplots}
\usetikzlibrary{calc,intersections}

\usetikzlibrary{shapes.geometric, arrows, positioning}

\tikzset{
	mybox/.style  = {draw, rectangle, minimum width=4cm, minimum height=0.8cm, text centered, text width=4.4cm,   
		font=\normalsize},
	box/.style  = {draw, rectangle, minimum width=2.0cm, minimum height=0.6cm, text centered, text width=3.0cm,   
		font=\normalsize},
	myarrow/.style = {line width=0.2pt, draw=black, -triangle 60, postaction={draw, line width=0.2pt, shorten >=10pt,-}}
}

\tikzstyle{arrow} = [->, >=stealth, -triangle 60]

\allowdisplaybreaks

\makeatletter
\newcommand{\leqnomode}{\tagsleft@true}
\newcommand{\reqnomode}{\tagsleft@false}
\makeatother

\begin{document}

\title{Linear Convergence of ISTA and FISTA\thanks{This work was supported by Grant No.YSBR-034 of CAS and Grant No.12288201 of NSFC.}}

\author[ ]{Bowen Li\qquad Bin Shi\thanks{Corresponding author, Email: \url{shibin@lsec.cc.ac.cn}} \qquad Ya-xiang Yuan} 

\affil[ ]{State Key Laboratory of Scientific and Engineering Computing, Academy of Mathematics and Systems Science, Chinese Academy of Sciences, Beijing 100190, China}

\affil[ ]{University of Chinese Academy of Sciences, Beijing 100049, China}

\date\today

\maketitle

\begin{abstract}
In this paper, we revisit the class of iterative shrinkage-thresholding algorithms (ISTA) for solving the linear inverse problem with sparse representation, which arises in signal and image processing. It is shown in the numerical experiment to deblur an image that the convergence behavior in the logarithmic-scale ordinate tends to be linear instead of logarithmic, approximating to be flat.  Making meticulous observations, we find that the previous assumption for the smooth part to be convex weakens the least-square model. Specifically, assuming the smooth part to be strongly convex is more reasonable for the least-square model, even though the image matrix is probably ill-conditioned. Furthermore, we improve the pivotal inequality tighter for composite optimization with the smooth part to be strongly convex instead of general convex, which is first found in~\citep{li2022proximal}. Based on this pivotal inequality, we generalize the linear convergence to composite optimization in both the objective value and the squared proximal subgradient norm. Meanwhile, we set a simple ill-conditioned matrix which is easy to compute the singular values instead of the original blur matrix. The new numerical experiment shows the proximal generalization of Nesterov's accelerated gradient descent (\texttt{NAG}) for the strongly convex function has a faster linear convergence rate than ISTA. Based on the tighter pivotal inequality, we also generalize the faster linear convergence rate to composite optimization, in both the objective value and the squared proximal subgradient norm, by taking advantage of the well-constructed Lyapunov function with a slight modification and the phase-space representation based on the high-resolution differential equation framework from the implicit-velocity scheme. 
\end{abstract}



\section{Introduction}
\label{sec: intro}

In scientific and engineering computation, the class of so-called inverse problems is to infer and calculate the causes according to the observations or effects, which is just the reverse of the forward inference from the causes to the effects. Specifically, the inverse problems involve determining the parameters of a system that cannot be observed directly. Since the eighties of the last century, it has become the fastest growing field, regardless of theory or algorithm, due to its wide range of applications, such as signal and image processing, statistical inference, geophysics, astrophysics, and so on~\citep{engl1996regularization}.

\begin{figure}[htpb!]
\centering
\includegraphics[scale=0.26]{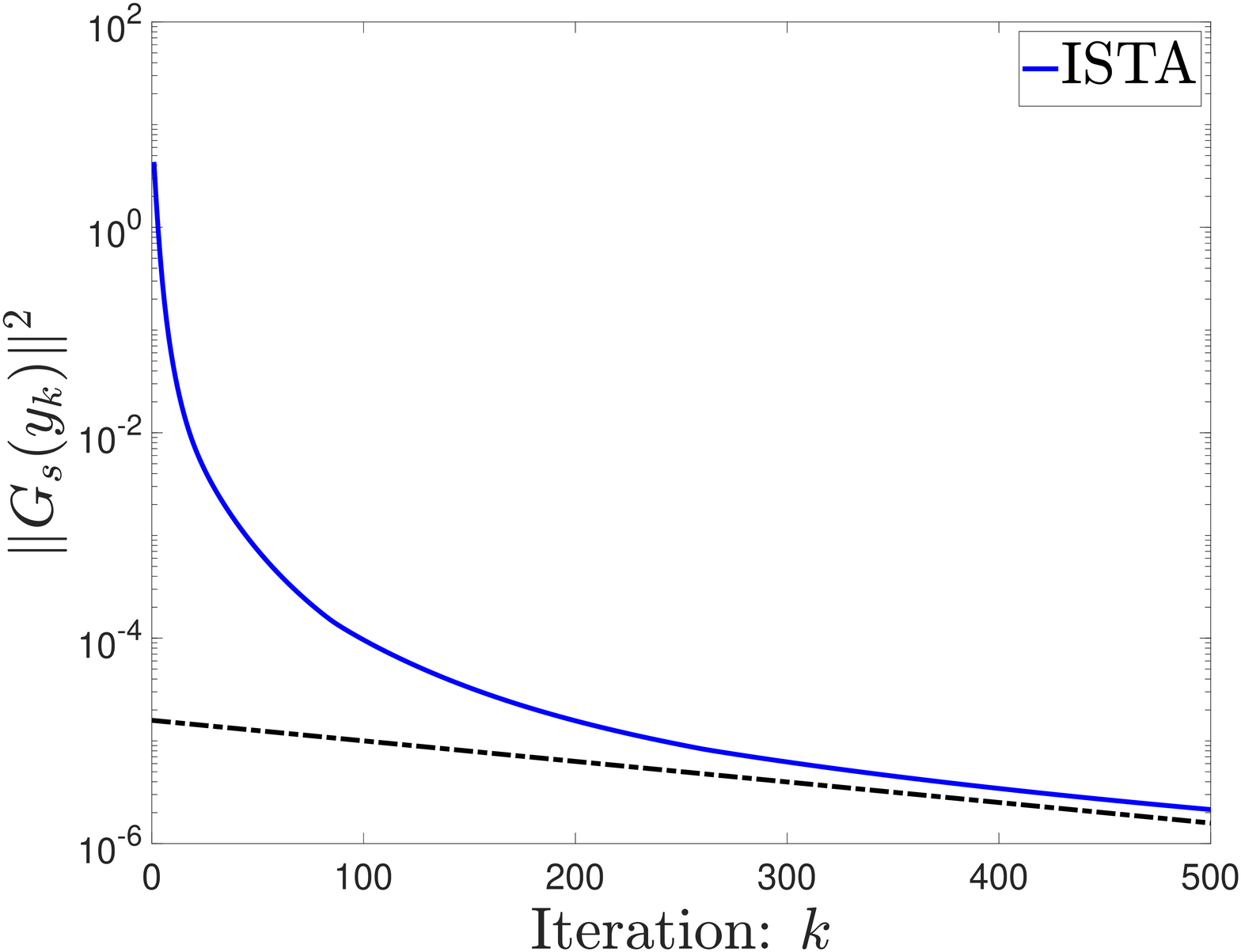}
\caption{The convergence behavior of the squared $s$-proximal subgradient norm by implementing ISTA to deblur the image of an elephant, shown in~\citep[Figure 1 \& Figure 2]{li2022proximal}.}
\label{fig: ista-elephant}
\end{figure}

Perhaps the simplest and classical one is the class of linear inverse problems, e.g., the image deblurring problem. Recall we deblur the image of an elephant with the pixel $256 \times 256$ in~\citep{li2022proximal}, which is implemented by the class of iterative shrinkage-thresholding algorithms (ISTA) as
\[
y_{k+1} = y_{k} - s G_s(y_{k}),
\]
with any initial $y_{0} = x_0 \in \mathbb{R}^d$ and a fixed step size $s \in (0,1/L]$, where $G_s(\cdot)$ is the $s$-proximal subgradient operator.\footnote{The Lipschitz constant $L$ and the $s$-proximal subgradient operator $G_s(\cdot)$ are defined rigorously in~\Cref{defn: L-smooth-strongly} and~\Cref{defn: proximal-subgradient}, respectively. }  We show the convergence behavior of the squared $s$-proximal subgradient norm in~\Cref{fig: ista-elephant}, where the regularization parameter is set as $\lambda = 10^{-6}$ and the step size $s = 0.5$. Based on the assumption of the convex function, we derive that the minimal squared $s$-proximal subgradient norm obeys the sublinear convergence, concretely $o(1/k^3)$,  in~\citep{li2022proximal}. However, with a black dashed line as a comparison, we find the convergence rate of the squared $s$-proximal subgradient norm tends to be linear instead of logarithmic $\sim -3\log k$, approximating to be flat. In other words, the objective function is more like the $\mu$-strongly convex function rather than the convex function. 

\paragraph{Numerical Phenomena} 
Recall the linear inverse problem with sparse representation is formulated by the least square (or linear regression) model with $\ell_1$-regularization in~\citep{beck2009fast} as
\begin{equation}
\label{eqn: lasso}
\frac{1}{2} \|Ax - b \|^2 +\lambda \|x\|_1,
\end{equation}
where $A \in \mathbb{R}^{m\times d}$ is an $m\times d$ matrix,  $b \in \mathbb{R}^m$ is an $m$-dimensional vector, and the regularization parameter $\lambda > 0$ is a tradeoff between fidelity to the measurements and noise sensitivity.\footnote{Throught the paper, $\|\cdot\|$ denotes the standard Euclidean norm without any special emphasization and labeling.  And  $\|\cdot\|_1$ denotes the $\ell_1$-norm as $\|x\|_1 = \sum_{i=1}^{d}|x_i|$ for any $x \in \mathbb{R}^d$. } Here for the least square model with $\ell_1$-regularization~\eqref{eqn: lasso}, we note $\lambda_i(A^{T}A)$, $(i=1,\ldots, d)$ as the eigenvalues for the symmetric matrix $A^TA$. Then, the parameters used for the optimization satisfy $\mu = \min_{1\leq i \leq d}\lambda_i(A^{T}A)$ and $L = \max_{1\leq i \leq d} \lambda_i(A^{T}A)$, respectively. In addition, the condition number is noted as
\[
\textbf{cond}(A^{T}A) = \frac{L}{\mu} = \frac{\max\limits_{1\leq i \leq d} \lambda_i(A^{T}A)}{\min\limits_{1\leq i \leq d}\lambda_i(A^{T}A)}.
\] 
Practically, the matrix $A$ in the least square model with $\ell_1$-regularization~\eqref{eqn: lasso} is ill-conditioned for the linear inverse problems~\citep{hansen2006deblurring, figueiredo2007gradient}, that is, the condition number $\textbf{cond}(A^{T}A)$ is very large.

\begin{figure}[htpb!]
\centering
\includegraphics[scale=0.26]{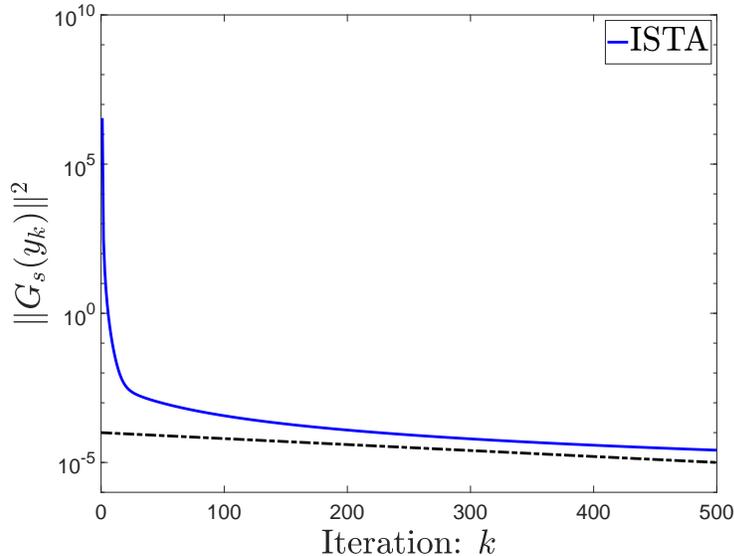}
\caption{The convergence behavior of the squared $s$-proximal subgradient norm by implementing ISTA in the least square model with $\ell_1$-regularization.}
\label{fig: ista}
\end{figure}

To further observe the phenomenon of the linear convergence discovered in~\Cref{fig: ista-elephant}, we conduct a new numerical experiment by setting the ill-conditioned matrix $A$ and the vector $b$ manually. Concretely, we set $A$ as a $500\times 500$ tridiagonal matrix with $A_{i,i} =2,\;A_{i,i+1}=1,\;A_{i+1,i}=1$ and $b$ as a $500$-dimensional vector with every entry $b_i=1$. Taking some simple computations, we obtain the optimization parameters as $\mu = \min_{1\leq i \leq d}\lambda_i(A^{T}A) = 1.5461 \times 10^{-9}$ and $L=\max_{1\leq i \leq d} \lambda_i(A^{T}A) = 15.9997$, respectively. Obviously, the matrix $A$ is ill-conditioned since the condition number is $\textbf{cond}(A^{T}A) \approx 1 \times 10^{-10}$. Meanwhile, this also manifests that setting the step size $s = 0.05 <1/L$ is reasonable. Finally, we show the numerical experiment by implementing ISTA with the same regularization parameter $\lambda = 10^{-6}$ in~\Cref{fig: ista}, where the linear convergence behavior of the squared proximal subgradient norm is more predominant after several iterations. 

\begin{figure}[htpb!]
\centering
\includegraphics[scale=0.26]{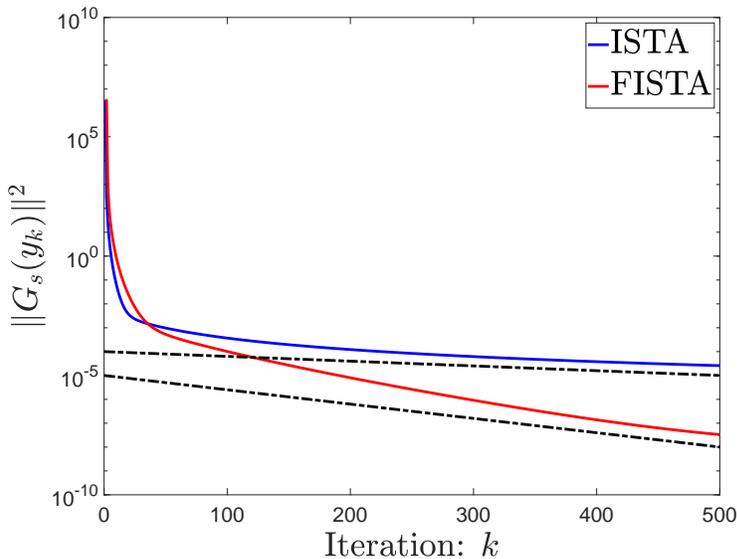}
\caption{The convergence behavior of the squared proximal subgradient norm by implementing ISTA and FISTA in the least square model with $\ell_1$-regularization.}
\label{fig: fista}
\end{figure}

For the smooth case, Nesterov's accelerated gradient descent~(\texttt{NAG}) globally accelerates the convergence rate of the vanilla gradient descent for the $\mu$-strongly convex function. Naturally, we hope to investigate the proximal version of~\texttt{NAG} as  
\[
\left\{\begin{aligned}
       & x_{k} =y_{k-1} - sG_s(y_{k-1}),              \\
       & y_{k} = x_{k}  + \frac{x_{k} - x_{k-1}}{1+2\sqrt{\mu s}},
       \end{aligned}\right.
\]
with any initial $x_0 = y_0 \in \mathbb{R}^d$, which is also called the class of fast iterative shrinkage-thresholding algorithms, shortened as FISTA.\footnote{The FISTA here is a proximal generalization of~\texttt{NAG} for the $\mu$-strongly convex function, and different from that in~\citep{beck2009fast}, which is based on~\texttt{NAG} for the convex function. In addition, the momentum coefficient here has a little modification, which is first proposed in~\citep{chen2022revisiting} to provide a tailor-made proof.} Thus, a new numerical experiment by adding FISTA for comparison is show in~\Cref{fig: fista}. Similarly, we also add a black dashed line as a basis to characterize the worst convergence rate of the squared $s$-proximal gradient norm by implementing FISTA. By comparing the two dashed lines in~\Cref{fig: fista}, we find that the acceleration phenomenon indeed exists for the squared $s$-proximal subgradient norm of FISTA beyond ISTA.

\subsection{Main contributions}
\label{subsec: overview-contribution}

From the description above, the linear inverse problem with sparse representation is formulated as the least-square model with $\ell_1$-regularization~\eqref{eqn: lasso}. In the language of optimization, it corresponds to composite optimization
\begin{equation}
\label{eqn: composite}
\min_{x\in\mathbb{R}^d}\Phi(x) := f(x) + g(x),
\end{equation}
where $f$ is a $\mu$-strongly smooth convex function and $g$ is a continuous convex function that is possibly not smooth. For the composite optimization problem~\eqref{eqn: composite}, we make the following contributions. 

\paragraph{The pivotal inequality for composite optimization} From the smooth case, the pivotal inequality for the convex function is observed  in~\citep{li2022proximal} as
\begin{equation}
\label{eqn: composite-smooth-c}
f(y - s \nabla f(y)) - f(x) \leq \langle \nabla f(y), y - x\rangle - \left(s- \frac{s^2L}{2}\right) \|\nabla f(y)\|^2.
\end{equation}
And then, we generalize~\eqref{eqn: composite-smooth-c} to composite optimization to derive the $s$-proximal subgradient norm minimization of ISTA and FISTA. Similarly, it is observed here that the pivotal inequality~\eqref{eqn: composite-smooth-c} can be improved tighter to the $\mu$-strongly convex function as
\begin{equation}
\label{eqn: composite-smooth-sc}
f(y - s \nabla f(y)) - f(x) \leq \langle \nabla f(y), y - x\rangle - \frac{\mu}{2} \| y - x\|^2 - \left(s- \frac{s^2L}{2}\right) \|\nabla f(y)\|^2. 
\end{equation}
We generalize the pivotal inequality~\eqref{eqn: composite-smooth-sc} to the composite objective function in~\Cref{sec: key-inequality}. 

\paragraph{Linear convergence of ISTA}
Recall the squared distance as the Lyapunov function in~\citep[Theorem A.4]{shi2019acceleration} is used to derive the linear convergence.\footnote{It should be noted that~\citet[Theorem 2.1.15]{nesterov1998introductory} use the squared norm  to derive the convergence rate, which is only recounted by the language of the Lyapunov function in~\citep[Theorem A.4]{shi2019acceleration}.} Based on the proximal version of the pivotal inequality~\eqref{eqn: composite-smooth-sc}, we obtain the convergence rate of ISTA for both the objective value and the $s$-proximal subgradient norm square, respectively as
\[
\Phi(y_{k}) - \Phi(x^\star)  \leq O\left(\left(\frac{1-\mu s}{1 + \mu s}\right)^k\right) \quad \text{and} \quad \|G_s(y_k)\|^2    \leq      O\left(\left(\frac{1-\mu s}{1 + \mu s}\right)^k\right)  
\] 
for any $0 < s \leq 1/L$. Moreover, the $s$-proximal subgradient norm is convenient to be used in practice, where the minimization rate is consistent with the numerical phenomena in~\Cref{fig: ista-elephant} and~\Cref{fig: ista}.


\paragraph{Accelerated linear convergence of FISTA}
Recall the Lyapunov function constructed from the implicit-velocity in~\citep{chen2022revisiting} is used to obtain the linear convergence. Here, we take a modification for the potential energy using $x_k$ instead of $y_k$. Combined with the proximal version of the pivotal inequality~\eqref{eqn: composite-smooth-sc}, we obtain the convergence rate of FISTA, respectively as:
\begin{itemize}
\item[\textbf{(1)}] For the objective value, it follows
\[
\Phi(x_{k}) - \Phi(x^\star)  \leq O\left( \frac{1}{\left(1 + \frac{\sqrt{\mu s}}{4}\right)^k  }\right) 
\]
with any step size $0 < s \leq 1/L$; 
\item[\textbf{(2)}] for the squared $s$-proximal subgradient norm, it follows 
\[
\|G_s(y_k)\|^2  \leq O\left( \frac{1}{\left(1 + \frac{\sqrt{\mu s}}{4}\right)^k  }\right) 
\]
with any step size $0<s<1/L$.
\end{itemize}
Similarly, the minimization rate of the squared $s$-proximal subgradient norm is consistent with the numerical phenomenon observed in~\Cref{fig: fista}.

%

\subsection{Related works and organization}
\label{subsec: notations-organization}

Recall the history of first-order smooth optimization,~\citet{nesterov1983method} proposed~\texttt{NAG} for the convex function and the ``\textit{Estimate Sequence}'' technique to bring about the global acceleration phenomenon, which was also generalized to the $\mu$-strongly convex function in~\citep{nesterov1998introductory}. Furthermore,~\citep{nesterov2010recent} points out that the algorithmic structures are very close in both~\texttt{NAG} for the $\mu$-strongly convex function and Polyak's heavy-ball method~\citep{polyak1964some}. And then, this view is emphasized further in~\citep{jordan2018dynamical}. Finally, the veil behind the acceleration phenomenon is lifted by the discovery of gradient correction based on the high-resolution differential equation framework in~\citep{shi2022understanding}. Recently, the high-resolution differential equation framework from the implicit-velocity scheme is found to be perfect and superior to the gradient-correction scheme in~\citep{chen2022revisiting}.

The technique of Lyapunov analysis used to investigate the acceleration phenomenon dates from~\citep{su2016differential} for the convex function. Besides capturing the acceleration mechanism, the Lyapunov function based on this high-resolution differential framework leads to the gradient norm minimization in~\citep{shi2022understanding}. The tailor-made proofs and the equivalent implicit-velocity effect are shown in~\citep{chen2022gradient}. With the observation to the key inequality, the gradient norm minimization is generalized to the composite function in~\citep{li2022proximal}. It is proposed by combining the low-resolution differential equation and the continuous limit of the Newton method to design and analyze new algorithms in~\citep{attouch2020first}, where the terminology is called inertial dynamics with a Hessian-driven term. Although the inertial dynamics with a
Hessian-driven term resembles closely with the high-resolution differential equations in~\citep{shi2022understanding}, it is important to note that the Hessian-driven terms are from the second-order information of Newton’s method~\citep{attouch2014dynamical}, while the gradient correction entirely relies on the first-order information of Nesterov's accelerated gradient descent method itself.

The remainder of this paper is organized as follows. In~\Cref{sec: prelim}, we collect some basic definitions and theorems as preliminaries. We provide two key inequalities for composite functions with rigorous proofs in~\Cref{sec: key-inequality}. The convergence rate of ISTA for both the objective value and the $s$-proximal subgradient norm square is shown in~\Cref{sec: nonaccele}. And we demonstrate the convergence rate of FISTA for both the objective value and the $s$-proximal subgradient norm square from the implicit-velocity scheme in~\Cref{sec: accele}. In~\Cref{sec: discussion}, we propose a conclusion and discuss some further research works.

\section{Preliminaries}
\label{sec: prelim}

In this section, we collect some basic definitions and theorems from the classical books~\citep{nesterov1998introductory}  and~\citep{rockafellar1970convex} with slight modifications for later uses in the following sections. First, some basic definitions in convex optimization are described from~\citep{nesterov1998introductory}, e.g., convex functions, subgradients, and subdifferentials.

\begin{defn}
\label{defn: convex}
A continuous function $g = g(x)$ defined on $\mathbb{R}^d$ is convex if for any $\alpha \in [0,1]$, the following inequlaity holds
\[
g\left( \alpha x + (1 - \alpha)y \right) \leq \alpha g(x) + (1 - \alpha)g(y)
\] 
for any $x, y \in \mathbb{R}^d$. Moreover, we note $\mathcal{F}^0(\mathbb{R}^d)$ as the class of continuous convex functions.
\end{defn}

\begin{defn}
\label{defn: subgradient}
Let $g \in \mathcal{F}^0(\mathbb{R}^d)$. A vector $v$ is called a subgradient of $g$ at $x$ if for any $y \in \mathbb{R}^d$, we have
\[
g(y) \geq g(x) + \langle v, y - x\rangle.
\] 
Moreover, the set of all subgradients of $g \in \mathcal{F}^0(\mathbb{R}^d)$ at $x$ is called the subdifferential of $g$ at $x$, and we note it as $\partial g(x)$. 
\end{defn}

Then, we collect two theorems from~\citep{rockafellar1970convex} to characterize the basic properties of subdifferentials and subgradients. 

\begin{theorem}[Theorem 23.8 in \citep{rockafellar1970convex}]
\label{thm: sum-diff}
Let $g_1, g_2 \in\mathcal{F}^0(\mathbb{R}^d)$. Then the subdifferentials satisfy
\[
\partial(g_1 + g_2)(x) = \partial g_1(x) + \partial g_2(x)
\]
for any $x \in \mathbb{R}^d$.
\end{theorem}

\begin{theorem}[Theorem 25.1 in \citep{rockafellar1970convex}]
\label{thm: subgrad-unique}
Let $g \in \mathcal{F}^0(\mathbb{R}^d)$. If $g$ is differentiable at $x$, then $\nabla g(x)$ is the unique subgradient of $g$ at $x$, so that in particular 
\[
g(y) \geq g(x) + \langle \nabla g(x), y - x\rangle,
\] 
for any $y \in \mathbb{R}^d$.
\end{theorem}

Next, we define the $\mu$-strongly convex function instead of the general convex function for the smooth part in the composite function to characterize the least-square model. 
\begin{defn}
\label{defn: L-smooth-strongly}
Let $f \in \mathcal{F}^0(\mathbb{R}^d)$. The function $f$ is $L$-smooth if $f$ is differentiable and the gradient of $f$ satisfies
\[
\|\nabla f(x) - \nabla f(y)\| \leq L \|x - y\|
\]
for any $x, y \in \mathbb{R}^d$. The function $f$ is $\mu$-strongly convex if $f$ is differentiable and  satisfies 
\[
f(y) \geq f(x) + \langle \nabla f(y), y - x \rangle + \frac{\mu}{2} \|y - x\|^2
\]
for any $x, y \in \mathbb{R}^d$. The subclass of $\mathcal{F}^0(\mathbb{R}^d)$ with the $L$-smoothness and $\mu$-strongly convexity is noted as $\mathcal{S}_{\mu,L}^1(\mathbb{R}^d)$.
\end{defn}


In this paper, we consider the composite function $\Phi= f + g$ with $f \in \mathcal{S}_{\mu,L}^1$ and $g \in \mathcal{F}^0$ and note $x^\star$ as its unique minimizer. Finally, we provide the definitions of both the $s$-proximal operator and the $s$-proximal subgradient operator originally shown in~\citep{beck2009fast, su2016differential} as 
\begin{defn}
\label{defn: proximal-subgradient}
For any $f\in \mathcal{S}_{\mu,L}^1$ and $g \in \mathcal{F}^0$, the $s$-proximal operator is defined as
\begin{equation}
\label{eqn: proximal-operator}
    P_s(x) := \mathop{\arg\min}_{y\in\mr^d}\left\{ \frac{1}{2s}\left\| y - \left(x - s\nabla f(x)\right) \right\|^2 + g(y) \right\},
\end{equation}
for any $x \in \mathbb{R}^d$. Furthermore, the $s$-proximal subgradient operator is defined as
\begin{equation}
\label{eqn: subgradient-operator}
G_s(x): = \frac{x - P_s(x)}{s}
\end{equation}
for any $x \in \mathbb{R}^d$.
\end{defn}
When $g$ is reduced to the $\ell_1$-norm, that is, $g(x) = \lambda \|x\|_1$,\ we can derive the closed-form expression of the $s$-proximal operator~\eqref{eqn: proximal-operator} at any $x \in \mathbb{R}^d$ for the special case as
\[
P_s(x)_i = \big(\left|\left(x - s\nabla f(x)\right)_i\right| - \lambda s\big)_+ \text{sgn}\big(\left(x - s\nabla f(x)\right)_i\big),
\]
where $i=1,\ldots,d$ is the index. 

%
\section{Two inequalities for composite optimization}
\label{sec: key-inequality}

In this section, based on the basic definitions and theorems shown in~\Cref{sec: prelim}, we generalize two basic inequalities of $f \in \mathcal{S}_{\mu,L}^1$ to the composite function $\Phi = f + g$ with $f \in \mathcal{S}_{\mu,L}^1$ and $g \in \mathcal{F}^0$.

\begin{lemma}
\label{lem: sc_proximal1}
Let $f \in \mathcal{S}_{\mu,L}^{1}(\mathbb{R}^d)$ and $g \in \mathcal{F}^0(\mathbb{R}^d)$. Then the composite function $\Phi = f + g$ satisfies
\begin{equation}
\label{eqn: sc_proximal1}
\Phi(x) - \Phi(x^\star) \geq \frac{\mu}{2} \|x - x^\star\|^2
\end{equation}
\end{lemma}

\begin{proof}[Proof of~\Cref{lem: sc_proximal1}]
With~\Cref{defn: convex} and~\Cref{defn: subgradient}, the basic convex inequality holds for any subgradient $v \in \partial\Phi(x^\star)$ as
\[
\Phi(x) \geq \Phi(x^\star) + \langle v, x - x^\star \rangle
\]
with any $x \in \mathbb{R}^d$. Since $x^\star$ is the unique global minimizer, we know that zero is a subgradient, that is, $0 \in \partial\Phi(x^\star)$. Meanwhile, due to $f \in \mathcal{S}_{\mu,L}^{1}(\mathbb{R}^d)$, we obtain the basic inequality for the $\mu$-strongly convex function as
\begin{equation}
\label{eqn:scv-1}
f(x) \geq f(x^\star) + \langle \nabla f(x^\star), x - x^\star \rangle + \frac{\mu}{2} \|x - x^\star\|^2
\end{equation}
with any $x\in \mathbb{R}^d$; and since $g \in \mathcal{F}^0(\mathbb{R}^d)$, the basic convex inequality holds for any subgradient $v'\in \partial g(x^\star)$ as
\begin{equation}
\label{eqn:cv-1}
g(x) \geq g(x^\star) + \langle v', x - x^\star \rangle
\end{equation}
with any $x \in \mathbb{R}^d$. Combined with~\eqref{eqn:scv-1} and~\eqref{eqn:cv-1}, we have
\begin{equation}
\label{eqn:scv-2}
\Phi(x) \geq \Phi(x^\star) + \langle \nabla f(x^\star) + v', x - x^\star \rangle + \frac{\mu}{2} \|x - x^\star\|^2
\end{equation}
with any $x \in \mathbb{R}^d$. Furthermore, with~\Cref{thm: sum-diff} and~\Cref{thm: subgrad-unique}, we obtain that for any $v \in \partial\Phi(x^\star)$, there exists a subgradient $v' \in \partial g(x^\star)$ such that $v = \nabla f(x^\star) + v'$.  With $0 \in \partial\Phi(x^\star)$, we complete the proof.
%
\end{proof}

Recall the proof of the pivotal inequality for composite optimization~\citep[Lemma 3.1]{li2022proximal}. Nevertheless, the function class here is $\mathcal{S}_{\mu,L}^{1}$ instead of $\mathcal{F}_L^1$, we can make the basic convex inequality in~\citep[(3.5a)]{li2022proximal} tighter as 
\[
f(x) \geq f(y) + \langle \nabla f(y), x - y \rangle + \frac{\mu}{2} \|x - y\|^2
\]
for any $x, y \in \mathbb{R}^d$, which brings about the pivotal inequality~\citep[(3.1)]{li2022proximal} to be tighter. Hence, we provide a rigorous description to conclude this section in the following theorem. 

\begin{lemma}
\label{lem: sc-proximal2}
Let $\Phi = f + g$ be a composite function with $f \in \mathcal{S}_{\mu,L}^{1}(\mathbb{R}^d)$ and $g \in \mathcal{F}^0(\mathbb{R}^d)$. Then, the following inequality
\begin{equation}
\label{eqn: sc-proximal2}
\Phi(y - G_s(y)) \leq \Phi(x) + \langle G_s(y), y - x\rangle - \left(s- \frac{s^2L}{2}\right) \|G_s(y)\|^2 - \frac{\mu}{2} \| y - x\|^2
\end{equation}
holds for any $x, y \in \mathbb{R}^d$.
\end{lemma}

\section{Linear convergence of ISTA}
\label{sec: nonaccele}

In this section, we investigate the convergence rate of ISTA for composite optimization in both the objective value and the squared $s$-proximal subgradient norm. Similar to the smooth case in~\citep[Theorem A.4]{shi2019acceleration}, we construct the Lyapunov function for the composite function as the squared distance, i.e., 
\begin{equation}
\label{eqn: ly-ista}
\mathcal{E}(k) = \|y_k - x^\star\|^2.\footnote{It should be noted that~\citet[Theorem 2.1.15]{nesterov1998introductory} use the squared norm  to derive the convergence rate, which is only recounted by the language of the Lyapunov function in~\citep[Theorem A.4]{shi2019acceleration}.}
\end{equation}
With the discrete Lyapunov function~\eqref{eqn: ly-ista}, we show the convergence rate in the following theorem.

\begin{theorem}
\label{thm: ista}
Let $\Phi = f+ g$ be a composite function with $f \in \mathcal{S}_{\mu,L}^{1}(\mathbb{R}^d)$ and $g \in \mathcal{F}^0(\mathbb{R}^d)$. Taking any step size $0 < s \leq 1/L$, the iterative sequence $\{y_k\}_{k=0}^{\infty}$ generated by ISTA satisfies
\begin{subequations}
  \begin{empheq}[left=\empheqlbrace]{align}
   \Phi(y_{k}) - \Phi(x^\star) & \leq \frac{1}{s}\left(\frac{1-\mu s}{1 + \mu s}\right)^k\|x_0 - x^\star\|^2, \label{eqn: ista_rate_obj}\\
    \|G_s(y_k)\|^2 &\leq \frac{4}{s^2}\left(\frac{1-\mu s}{1 + \mu s}\right)^k\|x_0 - x^\star\|^2. \label{eqn: ista_rate_gn}
  \end{empheq}
\end{subequations}
\end{theorem}

\begin{proof}[Proof of~\Cref{thm: ista}]
Putting the sequence $\{y_k\}_{k=0}^{\infty}$ generated by ISTA and the unique minimizer $x^\star$ into~\eqref{eqn: sc-proximal2}, ~\Cref{lem: sc_proximal1} brings about the following inequality as
\begin{equation}
\label{eqn: ista-1}
\Phi(y_{k+1}) - \Phi(x^\star) \leq \langle G_s(y_k), y_k - x^\star \rangle - \frac{s}{2} \|G_s(y_k)\|^2 - \frac{\mu}{2}\|y_k - x^\star\|^2.
\end{equation}
Furthermore, with the pivotal inequality~\eqref{eqn: ista-1}, we derive that the objective value satisfies the following inequality as 
\begin{equation}
\label{eqn: ista-2}
\Phi(y_{k+1}) - \Phi(x^\star)  \leq \frac{\|y_k - x^\star\|^2}{2s} - \frac{\|y_k - x^\star - s G_s(y_k)\|^2}{2s};
\end{equation}
and with the Cauchy-Schwarz inequality,  the squared $s$-proximal subgradient norm satisfies the following inequality as
\begin{equation}
\label{eqn: ista-3}
\|G_s(y_k)\|^2 \leq \frac{4\|y_k - x^\star\|^2}{s^2}.
\end{equation}
Then, we calculate the iterative difference of the Lyapunov function as
\begin{align*}
\mathcal{E}(k+1) - \mathcal{E}(k)& = \|y_{k+1} - x^\star\|^2 - \|y_{k} - x^\star\|^2 \\
                                 & = \|y_{k} - s G_s(y_k) - x^\star\|^2 - \|y_{k} - x^\star\|^2 \\
                                 & = - 2s \langle G_s(y_k), y_k - x^\star \rangle + s^2 \|G_s(y_k)\|^2 \\
                                 & \leq - 2s \left( \Phi(y_{k+1}) - \Phi(x^\star) \right) - s\mu \|y_k - x^\star\|^2
\end{align*}
where the last inequality follows from~\eqref{eqn: ista-1}. With~\Cref{lem: sc_proximal1}, we can obtain further
\[
\|y_{k+1} - x^\star\|^2 \leq (1 - s\mu)\|y_{k} - x^\star\|^2 - \mu s \|y_{k+1} - x^\star\|^2,
\]
which directly leads to the linear convergence of the Lyapunov function as 
\begin{equation}
\label{eqn: ista-4}
\mathcal{E}(k) = \|y_{k} - x^\star\|^2 \leq \left(\frac{1-\mu s}{1 + \mu s}\right)\|x_0 - x^\star\|^2.
\end{equation}
From~\eqref{eqn: ista-2},~\eqref{eqn: ista-3} and~\eqref{eqn: ista-4}, we derive the linear rates in~\eqref{eqn: ista_rate_obj} and~\eqref{eqn: ista_rate_gn}, conseqently the proof is complete.
\end{proof}

We conclude this section with the following corollary by setting the step size $s=1/L$.

\begin{corollary}
\label{coro: ista}
When the step size is set $s = 1/L$, the iterative sequence $\{y_k\}_{k=0}^{\infty}$ generated by ISTA satisfies
\begin{align}
    \Phi(y_{k}) - \Phi(x^\star)  &\leq L\left(\frac{L-\mu }{L + \mu }\right)^k\|x_0 - x^\star\|^2, \label{eqn: ista_rate_obj1}\\
    \|G_s(y_k)\|^2 &\leq 4L^2\left(\frac{L-\mu }{L + \mu }\right)^k\|x_0 - x^\star\|^2. \label{eqn: ista_rate_gn1}
\end{align}   
\end{corollary}
%
%
%
%
%
%
%

\section{Accelerated linear convergence of FISTA}
\label{sec: accele}
In this section, we investigate the convergence rate of FISTA for composite optimization in both the objective value and the squared $s$-proximal subgradient norm. Here, we choose the phase-space representation of FISTA from the implicit-velocity scheme as
\begin{equation}
\label{eqn: iv-phase}
\left\{
\begin{aligned}
           & x_{k+1} - x_{k} = \sqrt{s}v_{k+1},  \\
           & v_{k+1} - v_{k} = - \frac{2\sqrt{\mu s}v_{k}}{1+2\sqrt{\mu s}} - \sqrt{s}G_s(y_k),
\end{aligned}
\right.
\end{equation}
with any initial $x_0\in \mathbb{R}^d$ and $v_0 = 0$, where the iterative sequence $\{y_k\}_{k=0}^{\infty}$ satisfies
\begin{equation}
\label{eqn: iv-y-x}
y_k = x_{k} + \frac{\sqrt{s}v_{k}}{1+2\sqrt{\mu s}}.
\end{equation}
Similar to~\citep{chen2022gradient}, we move one space backward for the discrete Lyapunov function constructed in~\citep{chen2022revisiting} as
\begin{equation}
\mathcal{E}(k) = \underbrace{\Phi(x_{k}) - \Phi(x^{\star})}_{\mathcal{E}_{\textbf{pot}}(k)} + \underbrace{\frac{\|v_{k}\|^2}{4(1+2\sqrt{\mu s})^2}}_{\mathcal{E}_{\textbf{kin}}(k)} +  \underbrace{\frac{\left\|v_{k} +2\sqrt{\mu}(x_{k} - x^\star)\right\|^2 }{4}}_{\mathcal{E}_{\textbf{mix}}(k)}.\label{eqn: lyapunov-iv}
\end{equation}
In~\eqref{eqn: lyapunov-iv}, we use the iterative sequence $\{x_k\}_{k=0}^{\infty}$ instead of $\{y_k\}_{k=0}^{\infty}$ for the potential energy, which is because the proof for composite optimization is essentially dependent on the pivotal inequality~\eqref{eqn: sc-proximal2}. The convergence rate is formulated rigorously in the following theorem.

\begin{theorem}
\label{thm: fista}
Let $\Phi = f+ g$ be a composite function with $f \in \mathcal{S}_{\mu,L}^{1}(\mathbb{R}^d)$ and $g \in \mathcal{F}^0(\mathbb{R}^d)$. Taking any step size $0 < s \leq 1/L$, the iterative sequence $\{x_k\}_{k=0}^{\infty}$ generated by FISTA satisfies
\begin{equation}
\label{eqn: fista_rate_obj_1}
\Phi(x_{k}) - \Phi(x^\star) \leq \frac{\Phi(x_0) - \Phi(x^\star) + \mu\|x_0 - x^\star\|^2}{\left(1 + \frac{\sqrt{\mu s}}{4}\right)^k};
\end{equation}
if the step size satisfies $0 < s < 1/L$, the iterative sequence $\{y_k\}_{k=0}^{\infty}$ satisfies
\begin{equation}
\label{eqn: fista_rate_gn_1}
\|G_s(y_k)\|^2 \leq \frac{2\left(\Phi(x_0) - \Phi(x^\star) + \mu\|x_0 - x^\star\|^2\right)}{s(1-sL)\left(1 + \frac{\sqrt{\mu s}}{4}\right)^k}.
\end{equation}
\end{theorem}

\begin{proof}[Proof of~\Cref{thm: fista}]
First, we present the calculation of the Lyapunov function's iterative difference in three steps, similarly with~\citep{chen2022revisiting}. 
\begin{itemize}
\item[\textbf{(1)}] First, we start from the mixed energy $\mathcal{E}_{\textbf{mix}}(k)$. The critical point to calculating the iterative difference lies in
\[
\left[ v_{k+1}  + 2\sqrt{\mu}(x_{k+1} - x^\star) \right] -  \left[ v_{k}  + 2\sqrt{\mu}(x_{k} - x^\star) \right] =  - \sqrt{s} G_s(y_{k}),
\]
which directly follows the implicit-velocity phase-space representation~\eqref{eqn: iv-phase} and~\eqref{eqn: iv-y-x}. Then, we calculate the iterative difference of the mixed energy as
\begin{align*}
\mathcal{E}_{\textbf{mix}}(k+1) - \mathcal{E}_{\textbf{mix}}(k) =  & \frac{1}{4}\left\| v_{k+1}  + 2\sqrt{\mu}(x_{k+1} - x^\star) \right\|^2 -  \frac{1}{4}\left\| v_{k}  + 2\sqrt{\mu}(x_{k} - x^\star) \right\|^2 \\
=    & \frac12 \left\langle  - \sqrt{s} G_s(y_{k}),  v_{k}  + 2\sqrt{\mu}(x_{k} - x^\star) \right\rangle + \frac{s}{4}\| G_s(y_{k})\|^2  \\
=    & \frac12 \left\langle  - \sqrt{s} G_s(y_{k}),   \frac{v_{k}}{1+2\sqrt{\mu s}}  + 2\sqrt{\mu}(y_{k} - x^\star) \right\rangle + \frac{s}{4}\|G(y_{k})\|^2,
\end{align*}
where the last inequality also follows from~\eqref{eqn: iv-y-x}. Furthermore, with the second scheme of~\eqref{eqn: iv-phase}, we obtain the iterative difference as
\begin{equation}
\label{eqn: iv-1}
\mathcal{E}_{\textbf{mix}}(k+1) - \mathcal{E}_{\textbf{mix}}(k) = -  \frac{\sqrt{s}}{2} \left\langle  G_s(y_{k}), v_{k+1}\right\rangle - \sqrt{\mu s}  \left\langle G_s(y_{k}), y_{k} - x^\star \right\rangle - \frac{s}{4}   \|G_s(y_{k})\|^2. 
\end{equation}

\item[\textbf{(2)}] Then, we consider the kinetic energy  $\mathcal{E}_{\textbf{kin}}(k)$.  The detailed calculation of iterative difference is shown as
\begin{align}
\mathcal{E}_{\textbf{kin}}(k+1) - \mathcal{E}_{\textbf{kin}}(k) & =  \frac{\|v_{k+1}\|^2}{4(1+2\sqrt{\mu s})^2} -  \frac{\|v_{k}\|^2}{4(1+2\sqrt{\mu s})^2} \nonumber \\
            & =   \frac{\|v_{k+1}\|^2}{4(1+2\sqrt{\mu s})^2} - \frac{\|v_{k+1} + \sqrt{s} G_s(y_{k})\|^2}{4}  \nonumber \\
            & =  - \frac{\sqrt{s}}{2} \langle G_s(y_{k}), v_{k+1} \rangle - \frac{\sqrt{\mu s} \|v_{k+1}\|^2}{(1 + 2\sqrt{\mu s})^2} - \frac{s}{4} \| G_s(y_{k})\|^2, \label{eqn: iv-2}
\end{align}
where the second equality follows from the second scheme of~\eqref{eqn: iv-phase}.

\item[\textbf{(3)}] In the last step, we turn to the potential energy $ \mathcal{E}_{\textbf{pot}}(k)$. Recall in~\citep[Section 4.2]{chen2022revisiting}, the calculation of iterative difference is based on the basic $L$-smooth inequality. Here for the composite optimization, we plug $x_{k}$ and $x_{k+1}$ generated by FISTA into~\eqref{eqn: sc-proximal2} to derive the inequality as 
\begin{equation}
\label{eqn: fista-key-inq}
\Phi(x_{k+1}) - \Phi(x_{k}) \leq \langle G_s(y_k), y_{k} - x_{k} \rangle - \left(s- \frac{s^2L}{2}\right) \|G_s(y_k)\|^2,
\end{equation}
which is the essential difference from the smooth case. Then with~\eqref{eqn: fista-key-inq}, we calculate the iterative difference as
\begin{align}
\mathcal{E}_{\textbf{pot}}(k+1) - \mathcal{E}_{\textbf{pot}}(k)  & = \Phi(x_{k+1}) - \Phi(x_{k}) \nonumber\\
              & \leq \langle G_s(y_{k}), y_{k} - x_{k} \rangle - \left(s- \frac{s^2L}{2}\right) \|G_s(y_k)\|^2 \nonumber \\
              &= \langle G_s(y_{k}), x_{k+1} - x_{k} \rangle + s \|G_s(y_k)\|^2- \left(s- \frac{s^2L}{2}\right)\|G_s(y_k)\|^2 \nonumber \\
              &= \sqrt{s}\langle G_s(y_{k}), v_{k+1}\rangle + \frac{s}{2} \|G_s(y_k)\|^2 -\frac{s (1-sL) }{2} \|G_s(y_k)\|^2 , \label{eqn: iv-3}
\end{align}
where the second equality follows from the gradient step of FISTA. 
\end{itemize}

Using the three iterative differences~\eqref{eqn: iv-1},~\eqref{eqn: iv-2} and~\eqref{eqn: iv-3}, we obtain the iterative difference of the Lyapunov function as
\begin{equation}
\label{eqn: iv-id-lyapunov1}
\mathcal{E}(k+1) - \mathcal{E}(k)   \leq  -  \sqrt{\mu s}  \left\langle G_s(y_{k}), y_{k} - x^\star \right\rangle - \frac{(\sqrt{\mu s} + \mu s) \|v_{k+1}\|^2}{(1+2\sqrt{\mu s})^2} -\frac{s (1-sL) }{2} \|G_s(y_k)\|^2. 
\end{equation}
Putting the sequence $x_{k+1}$ generated by FISTA and the unique minimizer $x^\star$ into~\eqref{eqn: sc-proximal2}, we can see that~\Cref{lem: sc_proximal1} brings about the following inequality as
\begin{equation}
\label{eqn: fista1}
\Phi(x_{k+1}) - \Phi(x^\star) \leq \langle G_s(y_k), y_{k} - x^\star \rangle - \left(s - \frac{sL^2}{2}\right) \|G_s(y_k)\|^2 - \frac{\mu}{2}\|y_k - x^\star\|^2.
\end{equation}
It follows from~\eqref{eqn: fista1} and~\eqref{eqn: iv-id-lyapunov1} that
\begin{multline}
\mathcal{E}(k+1) - \mathcal{E}(k)    \\
                                                      \leq -\sqrt{\mu s} \left[ \Phi(x_{k+1}) - \Phi(x^\star) + \frac{(1 + \sqrt{\mu s} ) \|v_{k+1}\|^2}{(1+2\sqrt{\mu s})^2}  + \frac{\mu \|y_{k} - x^\star\|^2 }{2} \right] -\frac{s (1-sL) }{2} \|G_s(y_k)\|^2. \label{eqn: iv-id-lyapunov2} 
\end{multline}
With the gradient step of FISTA, we estimate the mixed energy by Cauchy-Schwarz inequality as
\begin{align*}
\|v_{k+1} +2\sqrt{\mu}(x_{k+1} - x^\star) \|^2  & =     \left\| v_{k+1} +2\sqrt{\mu}(y_k - sG_s(y_k) - x^\star)  \right\|^2 \\
                                                & \leq  2 \|v_{k+1}\|^2 + 8 \mu \|y_k - sG_s(y_k) - x^\star\|^2             \\ 
                                                & =     2 \|v_{k+1}\|^2 + 8 \mu \|y_k  - x^\star\|^2 +  8 \mu s^2 \| G_s(y_k) \|^2 - 16 \mu s \langle G_s(y_k), y_k - x^\star \rangle \\
                                                & \leq  2 \|v_{k+1}\|^2 + 8 \mu \|y_k  - x^\star\|^2,     
\end{align*}
where the last inequality follows~\eqref{eqn: fista1} and~\eqref{eqn: sc_proximal1}. Furthermore, we estimate the discrete Lyapunov function~\eqref{eqn: lyapunov-iv} as
\begin{equation}
\label{eqn: lyapunov-iv-estimate}
\mathcal{E}(k)  \leq \Phi(x_{k}) - \Phi(x^{\star}) +   \frac{\left(\frac{3}{4}+\sqrt{\mu s} + \mu s\right)}{(1+2\sqrt{\mu s})^2} \|v_{k}\|^2  + 2\mu \|y_{k-1} - x^\star\|^2. 
\end{equation}
Comparing the coefficients of the corresponding terms in~\eqref{eqn: iv-id-lyapunov2} and~\eqref{eqn: lyapunov-iv-estimate} for $\mathcal{E}(k+1)$, we
obtain that the iterative difference of the discrete Lyapunov function satisfies
\[
\mathcal{E}(k+1) - \mathcal{E}(k) \leq -\sqrt{\mu s} \min\left\{1, \frac{1+\sqrt{\mu s}}{\frac{3}{4} + \sqrt{\mu s} + \mu s}, \frac14 \right\}\mathcal{E}(k+1) -\frac{s (1-sL) }{2} \|G_s(y_k)\|^2.
\]
We conclude here with the initial values $x_0 \in \mathbb{R}^d$ and $v_0 = 0$. 
\end{proof}

If we retrogress to the smooth case, the Lyapunov function~\eqref{eqn: lyapunov-iv} constructed by using $x_k$ instead of $y_k$ becomes more perfect than~\citep[Theorem 4.2]{chen2022revisiting} in terms of the initial values.  Taking a comparison between~\Cref{thm: ista} and~\Cref{thm: fista}, we can find a little difference for the constant in the convergence rate. However, different from the smooth case, we cannot use $L\|x_0 - x^\star\|^2$ to bound $\Phi(x_0) - \Phi(x^\star) + \mu \|x_0 - x^\star\|^2$ obtained in~\eqref{eqn: fista_rate_obj_1} and~\eqref{eqn: fista_rate_gn_1} for the composite function, since the $L$-smooth inequality does not work anymore.  To make the constant in the convergence rate consistent with~\Cref{thm: ista}, we need to start from $\mathcal{E}(1)$ in~\eqref{eqn: lyapunov-iv} as
\begin{align}
\mathcal{E}(1) & =    \Phi(x_1) - \Phi(x^\star) + \frac{\|v_1\|^2}{4(1+2\sqrt{\mu s})^2} + \frac{\|v_1 + 2\sqrt{\mu}(x_1 - x^\star)\|^2}{4}  \nonumber \\
               & \leq \Phi(x_1) - \Phi(x^\star) +  \frac{\left(\frac{3}{4}+\sqrt{\mu s} + \mu s\right)}{(1+2\sqrt{\mu s})^2}\|v_1\|^2 + 2\mu\|x_0 - x^\star\|^2, \label{eqn: initial-fista-1}
\end{align}
where the inequality follows from~\eqref{eqn: lyapunov-iv-estimate}. Taking some basic calculations, we obtain $v_1 = \sqrt{s}G_s(x_0)$ and put it into~\eqref{eqn: initial-fista-1}  as 
\begin{equation}
\label{eqn: initial-fista-2}
\mathcal{E}(1) \leq \Phi(x_1) - \Phi(x^\star) +  \frac{\left(\frac{3}{4}+\sqrt{\mu s} + \mu s\right)}{(1+2\sqrt{\mu s})^2}\cdot s\|G_s(x_0)\|^2 + 2\mu\|x_0 - x^\star\|^2.
\end{equation}
Similarly, putting  $x_1$ generated by FISTA and the unique minimizer $x^\star$ into~\eqref{eqn: sc-proximal2}, we can see that~\Cref{lem: sc_proximal1} brings about the following inequality as
\[
\Phi(x_1) - \Phi(x^\star) \leq \langle G_s(x_0), x_0 - x^\star \rangle - \frac{s}{2} \|G_s(x_0)\|^2.
\]
Furthermore, with perfect square and Cauchy-Schwarz inequality, we have the following two inequalities 
\[
\Phi(x_1) - \Phi(x^\star)  \leq \frac{\|x_0 - x^\star\|^2}{2s}  \qquad \text{and} \qquad \|G_s(x_0)\|^2  \leq \frac{4 \|x_0 - x^\star\|^2}{s^2}
\]
holds for any $0 < s \leq 1/L$. Hence, we conclude this section with the following corollary.

\begin{corollary}
\label{coro: fista1}
Let $\Phi = f+ g$ be a composite function with $f \in \mathcal{S}_{\mu,L}^{1}(\mathbb{R}^d)$ and $g \in \mathcal{F}^0(\mathbb{R}^d)$. If the step size $0 < s \leq 1/L$, the iterative sequence $\{y_k\}_{k=0}^{\infty}$ generated by FISTA satisfies
\begin{equation}
\label{eqn: fista_obj_2}
\Phi(x_{k}) - \Phi(x^\star) \leq \frac{11\|x_0 - x^\star\|^2}{2s\left(1 + \frac{\sqrt{\mu s}}{4}\right)^k};
\end{equation}
if the step size is set $s=1/L$, we have
\begin{equation}
\label{eqn: fista_obj_2l}
\Phi(x_{k}) - \Phi(x^\star) \leq \frac{11L\|x_0 - x^\star\|^2}{2\left(1 + \frac{1}{4}\sqrt{\frac{\mu}{L}}\right)^k}.
\end{equation}
Furthermore, if $0 < s < 1/L$, the iterative sequence $\{y_k\}_{k=0}^{\infty}$ satisfies
\begin{equation}
\label{eqn: fista_rate_gn_2}
\|G_s(y_k)\| \leq \frac{11\|x_0 - x^\star\|^2}{s^2(1-sL)\left(1 + \frac{\sqrt{\mu s}}{4}\right)^k}.
\end{equation}
\end{corollary}

\section{Conclusion and discussion}
\label{sec: discussion}

\begin{figure}[htpb!]
\centering
\includegraphics[scale=0.26]{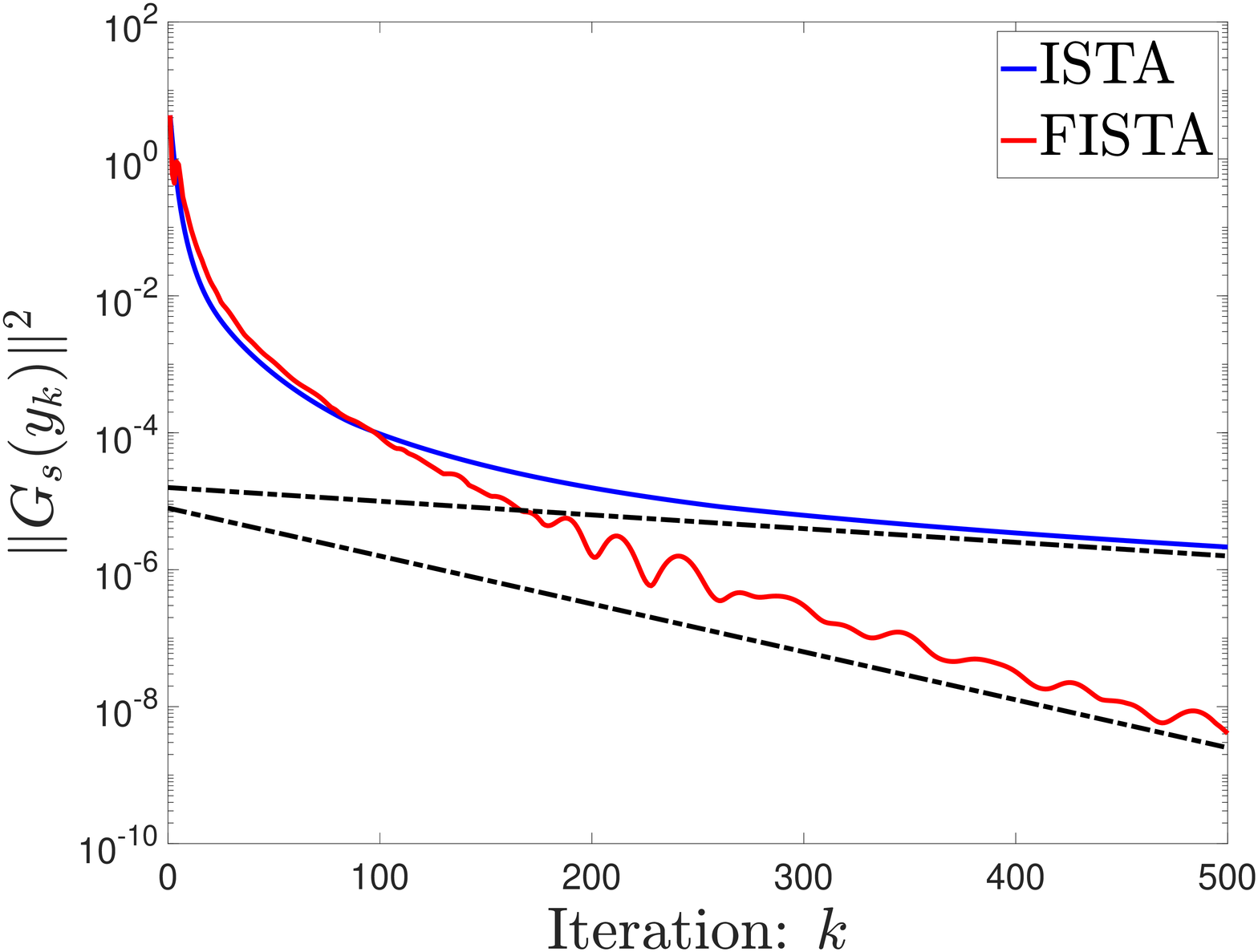}
\caption{The convergence behavior of the squared $s$-proximal subgradient norm by implementing ISTA and FISTA to deblur the image of an elephant, previously shown in~\citep[Figure 1 \& Figure 2]{li2022proximal}.}
\label{fig: fista1}
\end{figure}

In this paper,  we improve the pivotal inequality for composite optimization found in~\citep[Lemma 3.1]{li2022proximal}, which becomes tighter with the convexity of the smooth part substituted by the $\mu$-strongly convexity. Then, we apply the Lyapunov analysis to obtain the linear convergence for ISTA in both the objective value and the squared $s$-proximal subgradient norm. Meanwhile,  we also use the Lyapunov analysis via the phase-space representation from the implicit-velocity scheme to derive the accelerated linear convergence for FISTA in both the objective value and the squared $s$-proximal subgradient norm.  Looking back to the numerical experiment of FISTA in~\Cref{fig: fista},  we can compute the constant $\mu$ to set the momentum coefficient.  However, for the image deblurring problem in practice, it is so hard to compute the constant $\mu$ directly that we can set the momentum coefficient only empirically.  Nevertheless,  when we set the constant $\mu = 0.001$ in the numerical experiment of deblurring the image of an elephant, FISTA indeed converges faster than ISTA, which is shown in~\Cref{fig: fista1}. Therefore, how to provide an available constant $\mu$ in practice, e.g., for the blur matrix, looks very interesting. 

The $l_1$-norm is widely used in the statistics, such as the Lasso model and its variants, as the regularizer to capture the sparse structure. An interesting direction is to investigate further the theory in statistics related to Lasso based on the new discoveries of the composite optimization theory. For example, the sorted $l_1$-penalized estimation (SLOPE) is proposed to capture adaptively unknown sparsity in the asymptotic minimax problem~\citep{su2016slope},  and the Lasso path involves a strong relation to false discoveries  in~\citep{su2017false}. In addition,   the differential inclusion used for sparse recovery~\citep{osher2016sparse} is also worth further consideration in the future since is strongly related to composite optimization.

{\small
\subsection*{Acknowledgments}
We would like to thank Shuo Chen for his helpful discussions.
\bibliographystyle{abbrvnat}
\bibliography{reference}
}
\end{document}